\documentclass[english]{kart}
\usepackage{babel}
\usepackage{epsfig,graphicx,times}   
\usepackage{amssymb}
\newcounter{theoremcounter}
\def\Ps{\mathcal{P}}
\newcommand{\T}{\hspace*{0.6cm}}
\newcommand{\car}{\mbox{char}}

\def\di{m}
\newcommand{\eop}{\hfill{$\Box$}}
\def\thetheoremcounter{\arabic{theoremcounter}}
\newcommand{\labell}[1]{\label{#1}%
	\ifmmode $$\vspace*{-\baselineskip}\marginpar{#1}%
	\vspace*{-\baselineskip}$$\else\marginpar{#1}\fi}
\newenvironment{satz}%
	{\begin{trivlist}\refstepcounter{theoremcounter}%
	\item[]{{\bf Theorem \thetheoremcounter\ }}\hspace{2mm}\em}%
	{\end{trivlist}}
\newenvironment{lemma}%
	{\begin{trivlist}\refstepcounter{theoremcounter}%
	\item[]{{\bf Lemma \thetheoremcounter\ }}\hspace{2mm}\em}%
	{\end{trivlist}}	 
\newenvironment{definition}%
	{\begin{trivlist}\refstepcounter{theoremcounter}%
	\item[]{{\bf Definition
	\thetheoremcounter
	\ }}\hspace{2mm}\em}%
	{\end{trivlist}}	 	 
\newenvironment{algorithmus}%
	{\begin{trivlist}\refstepcounter{theoremcounter}%
	\item[]{\bf Algorithm \thetheoremcounter\ }\hspace{2mm}}%
	{\end{trivlist}}	 
\newenvironment{proof}%
	{\begin{trivlist}
	\item[]{{\sc Proof.}}
	}%
	{\eop\noindent\end{trivlist}}	 
\begin{document}
\begin{frontmatter}
\title{On the characteristic of integral point sets in $\mathbb{E}^\di$}
\author{Sascha Kurz}
\ead{sascha.kurz@uni-bayreuth.de}
\ead[url]{www.wm.uni-bayreuth.de}
\address{University of Bayreuth, Department of Mathematics, D-95440 Bayreuth, Germany}
\begin{abstract}
	We generalise the definition of the characteristic of an integral triangle to integral simplices and 
	prove that each simplex in an integral point set has the same characteristic. This theorem is used for 
	an efficient construction algorithm for integral point sets. Using this algorithm we are able to 
	provide new exact values for the minimum diameter of integral point sets.
\end{abstract}
\end{frontmatter}
	\parskip1.3em
	\section{Introduction}
	Since the time of the Pythagoreans, mathematicians have considered geometrical objects with integral sides.  
	Here we study sets of points in the Euclidean space $\mathbb{E}^\di$ where the pairwise distances are
	integers. Although there is a long history for integral point sets, very little is known about integral 
	point sets for dimension $m\ge 3$, see \cite{integral_distances_in_point_sets} for an overview.

	Due to Heron  the area of a triangle with side lengths $a$, $b$, and $c$ is given 
	by $$A_\Delta=\frac{\sqrt{(a+b+c)(a+b-c)(a-b+c)(-a+b+c)}}{4}.$$ Thus we can write 
	the area as $A_\Delta=q\sqrt{k}$ with a rational number $q$ and a squarefree integer 
	$k$. If $A_\Delta\neq 0$ the integer $k$ is unique and is called the \textbf{characteristic} 
	or the \textbf{index} of the triangle. This invariant receives its relevance from the 
	following theorem \cite{hab_kemnitz}.
  
	\begin{satz}{}
		\label{satz_char_2}
		The triangles spanned by each three non collinear points in a plane integral point set 
		have the same characteristic.	
	\end{satz}
  
	This theorem can be utilised to develop an efficient algorithm for the generation of plane integral 
	point sets \cite{phd_kurz,kurz_wassermann}. Here we will generalise the definition of the characteristic of an integral 
	triangle to integral simplices and prove an analogue to Theorem \ref{satz_char_2}.
	Later on we will use this theorem to develop a generation algorithm for integral point 
	sets in $\mathbb{E}^\di$ and present some new numerical data.
  
	
	\section{Characteristic of integral simplices}
	As the definition of the characteristic of an integral triangle depends on the 
	area of a triangle we consider the volume of an $\di$-dimensional
	simplex for point sets  in $\mathbb{E}^\di$. Therefore we need the Cayley-Menger matrix
	of a point set.
  
	\begin{definition}
		If $\Ps$ is a point set in $\mathbb{E}^\di$ with vertices $v_0,v_1,\dots,v_{n-1}$ and 
		$C=(d_{i,j}^2)$ denotes the $n\times n$ matrix given by $d_{i,j}^2=\Vert v_i-v_j\Vert_2^2$ 
		the Cayley-Menger matrix $\hat{C}$ is obtained from $C$ by bordering $C$ with a top 
		row $(0,1,1,\dots,1)$ and a left column $(0,1,1,\dots,1)^T$.
  \end{definition}
  
	By $CMD(\{v_0,v_1,\dots,v_{n-1}\})$ we denote the determinant of $\hat{C}(\{v_0,v_1,\dots,v_{n-1}\})$.
	If $n=\di+1$, the $\di$-dimensional volume $V_{\di}$ of $\Ps$ is given by
	$$
		V_{\di}(\Ps)^2=\frac{(-1)^{\di+1}}{2^\di(\di!)^2}\det(\hat{C})\,.
	$$
	This allows us to define the characteristic of an $\di$-dimensional integral simplex 
	to be the squarefree integer $k$ in $V_{\di}(\Ps)=q\sqrt{k}$ whenever 
	$V_{\di}(\Ps)\neq 0$ and $q\in\mathbb{Q}$. In order to prove the proposed theorem 
	we consider a special coordinate representation of integral simplices. 
    
	\begin{lemma}
		\label{lemma_coordinates}
		An integral $\di$-dimensional simplex $\mathcal{S}=\{v_0',v_1',\dots,v_{\di}'\}$ with 
		distance matrix $D=(d_{i,j})_{0\le i,j\le\di}$ and $V_{\di}(\mathcal{S})\neq 0$ can be 
		transformed via an isometry into the coordinates
		\begin{eqnarray*}
			v_0&=&(0,0,\dots,0),\\
			v_1&=&(q_{1,1}\sqrt{k_1},0,0\dots,0),\\
			v_2&=&(q_{2,1}\sqrt{k_1},q_{2,2}\sqrt{k_2},0,\dots,0),\\
			\vdots\\
			v_{\di}&=&(q_{\di,1}\sqrt{k_1},q_{\di,2}\sqrt{k_2},\dots,q_{\di,\di}\sqrt{k_{\di}}),
		\end{eqnarray*}
		where $k_i$ is the squarefree part of 
		$\frac{V_i(v_0',v_1',\dots,v_i')^2}{V_{i-1}(v_0',v_1',\dots,v_{i-1}')^2}$, 
		$q_{i,j}\in\mathbb{Q}$, and $q_{j,j},k_j\neq 0$.
  \end{lemma}
  \begin{proof}
  		We can obviously set $v_0=(0,0,\dots,0)$ and since $d_{0,1}\in\mathbb{N}$ 
		we can furthermore set $v_1=(d_{0,1}\sqrt{k_1},0,0,\dots,0)$ where 
		$k_1=\frac{V_1(v_0',v_1')}{V_0(v_0')}=1$. Now we assume that 
		we have already transformed $v_0',v_1',\dots,v_{i-1}'$ into the stated coordinates.
		We set $v_i=(x_1,x_2,\dots,x_{\di})$ with $x_j\in\mathbb{R}$. Since the points 
		$v_0,v_1,\dots,v_i$ span an $i$-dimensional hyperplane of $\mathbb{E}^{\di}$ we 
		can set $x_{i+1}=\dots=x_{\di}=0$. For $j\le i$ we have 
		$$
			d_{j,i}^2=\Vert v_j-v_i\Vert_2^2=\sum_{h=1}^j(q_{j,h}\sqrt{k_h}-x_h)^2+
			\sum_{h=j+1}^ix_h^2\,.
		$$
		For $0< j<i$ we consider
		$$
			d_{0,i}^2-d_{j,i}^2=\sum_{h=1}^j x_h^2- (q_{j,h}\sqrt{k_h}-x_h)^2
		$$ 
		where we can set $x_h=q_{i,h}\sqrt{k_h}$ for $h<j$ by induction, yielding
		$$
			d_{0,i}^2-d_{j,i}^2=-q_{j,j}^2k_h+2q_{j,j}\sqrt{k_h}x_j+\sum_{h=1}^{j-1} 2q_{i,h}q_{j,h}k_h-q_{j,h}^2k_h\,.
		$$
		Thus 
		$$
			x_j=\frac{q_{j,j}^2k_h+\sum\limits_{h=1}^{j-1}(q_{j,h}^2k_h-2q_{i,h}q_{j,h}k_h)+d_{0,i}^2-d_{j,i}^2}{2q_{j,j}\sqrt{k_h}}
		$$
		and we can write $x_j=q_{i,j}\sqrt{k_j}$ since $2q_{j,j}\sqrt{k_h}\neq 0$ due to induction. 
		With this we have 
		$$
			d_{0,i}^2=\sum_{h=1}^ix_h^2=x_i^2+\sum_{h=1}^{i-1}q_{i,h}^2k_h.
		$$
		Thus 
		$$
			x_i=\sqrt{d_{0,i}^2-\sum_{h=1}^{i-1}q_{i,h}^2k_h}=q_{i,i}\sqrt{k_i}\,.
		$$
		We also have $q_{i,i}\sqrt{k_i}\neq 0$ since $v_0',v_1',\dots,v_i'$ cannot lie 
		in an $i-1$-dimensional hyperplane of $\mathbb{E}^{\di}$ due to 
		$V_{\di}(v_0',v_1',\dots,v_{\di}')\neq 0$.
	\end{proof}
  
	The $k_j$ are associated to the characteristic $\car(\mathcal{S})=k$ in the following way
	$$
		\car(\mathcal{S})=k=\mbox{squarefree part of }\prod_{j=1}^{\di}k_j\,.
	$$
  
	\begin{satz}{}
		\label{characteristic}
		In an $\di$-dimensional integral point set $\Ps$ all simplices $\mathcal{S}=\{v_0,v_1,\dots,v_{\di}\}$
		with $V_{\di}(\mathcal{S})\neq 0$ have the same characteristic $\car(S)=k$.
	\end{satz}
	\begin{proof}
		It suffices to prove that $\car(\mathcal{S}_1)=\car(\mathcal{S}_2)$ for two integral simplices
		$\mathcal{S}_1=\{v_0,v_1,\dots.v_{\di}\}$ and $\mathcal{S}_2=\{v_0,\dots,v_{\di-1},v_{\di}'\}$ 
		with $V_{\di}(\mathcal{S}_1),V_{\di}(\mathcal{S}_2)\neq 0$. With the notations from Lemma 
		\ref{lemma_coordinates} we have  for the distance between $v_{\di}$ and $v_{\di}'$,
		\begin{eqnarray*}
			d(v_{\di},v_{\di}')^2&=&\sum_{i=1}^{\di}(q_{\di,i}\sqrt{k_i}-q_{\di,i}'\sqrt{k_i'})^2\\
			&=&\sum_{i=1}^{\di}(q_{\di,i}\sqrt{k_i}-q_{\di,i}'\sqrt{k_i})^2+(q_{\di,\di}\sqrt{k_\di}-q_{\di,\di}'
			\sqrt{k_\di'})^2\\ 
			&=&\sum_{i=1}^{\di-1}(q_{\di,i}-q_{\di,i}')^2k_i+q_{\di,\di}^2k_{\di}-2q_{\di,\di}q_{\di,\di}'
			\sqrt{k_{\di}k_{\di}'}+q_{\di,\di}'^2k_{\di}'\,.
		\end{eqnarray*}
		Thus $\sqrt{k_{\di},k_{\di}'}$ has to be an integer. Because $k_{\di}$ and $k_{\di}'$ are squarefree integers 
		$\neq 0$ we have $k_{\di}=k_{\di}'$ and so $\car(\mathcal{S}_1)=\car(\mathcal{S}_2)$.
	\end{proof}
  
  
	\section{Construction of integral point sets}
  
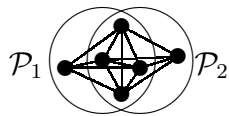
\begin{figure}[h]
		\begin{center}
			\setlength{\unitlength}{0.5cm}
			\begin{picture}(5,4)
				\put(2,2){\circle{4}}
				\put(3,2){\circle{4}}
				\put(2,2){\circle*{0.4}}
				\put(2.5,1.1){\circle*{0.4}}
				\put(2.5,2.9){\circle*{0.4}}
				\put(3,1.8){\circle*{0.4}}
				\put(1,1.8){\circle*{0.4}}
				\put(4,2.1){\circle*{0.4}}
				\put(0,0){\qbezier(2,2)(2.25,1.55)(2.5,1.1)}
				\put(0,0){\qbezier(2,2)(2.25,2.45)(2.5,2.9)}
				\put(0,0){\qbezier(2,2)(2.5,1.9)(3,1.8)}
				\put(0,0){\qbezier(2,2)(1.5,1.9)(1,1.8)}
				\put(0,0){\qbezier(2.5,1.1)(2.5,2.0)(2.5,2.9)}
				\put(0,0){\qbezier(2.5,1.1)(2.75,1.45)(3,1.8)}
				\put(0,0){\qbezier(2.5,1.1)(1.75,1.45)(1,1.8)}
				\put(0,0){\qbezier(2.5,2.9)(2.75,2.35)(3,1.8)}
				\put(0,0){\qbezier(2.5,2.9)(1.75,2.35)(1,1.8)}
				\put(0,0){\qbezier(3,1.8)(2,1.8)(1,1.8)}
				\put(0,0){\qbezier(4,2.1)(3,2.05)(2,2)}
				\put(0,0){\qbezier(4,2.1)(3.25,1.6)(2.5,1.1)}
				\put(0,0){\qbezier(4,2.1)(3.25,2.5)(2.5,2.9)}
				\put(0,0){\qbezier(4,2.1)(3.5,1.95)(3,1.8)}
				\put(-0.5,1.65){$\mathcal{P}_1$}
				\put(4.45,1.65){$\mathcal{P}_2$}
			\end{picture}
		\end{center}
		\caption{Combination of two integral point sets.}
		\label{fig2}
	\end{figure}

	The key principle for a recursive construction of integral point set consisting of $n$ points is the 
	combination of two integral point sets $\mathcal{P}_1=\{v_0,\dots,v_{n-2}\}$ and $\mathcal{P}_2=\{v_0,\dots,v_{n-3},v_{n-1}\}$ 
	consisting of $n-1$ points sharing $n-2$ points, see Figure \ref{fig2}. Here we describe an integral point set by 
	a symmetric matrix $D=(d_{i,j})$ 
    
	representing the distances between the points. Because not all symmetric 
	matrices  are realizable as distance matrices in $\mathbb{E}^{\di}$ we need a generalisation of the triangle 
	inequalities.
  
	\begin{satz}{\textbf{(Menger \cite{menger})}}
		\label{theorem_menger}
		A set of vertices $\{v_0,v_1,\dots,v_{n-1}\}$ with pairwise distances $d_{i,j}$ is realizable 
		in the Euclidean space $\mathbb{E}^{\di}$ if and only if for all subsets $\{i_0,i_1,\dots,i_{r-1}\}
		\subset\{0,1,\dots,n-1\}$ of cardinality $r\le \di+1$,
		$$
			(-1)^{r}CMD(\{v_{i_0},v_{i_1},\dots,v_{i_{r-1}}\})\ge 0,
		$$
		and for all subsets of cardinality $\di+2\le r\le n$,
		$$
			(-1)^{r}CMD(\{v_{i_0},v_{i_1},\dots,v_{i_{r-1}}\})=0\,.
		$$
	\end{satz}
  
	Fortunately we do not need to check all these equalities and inequalities. Because the point sets
	$\Ps_1$ and $\Ps_2$ are realizable due to our construction strategy it suffices to check  
	$(-1)^{n}CMD(\{v_0,v_1,\dots,v_{n-1}\})$ \cite{phd_kurz}. 
  
	To solve the equivalence problem for integral point sets we use  a variant of orderly 
	generation \cite{0412.05006,orderly,mckay,winner}. For the required ordering we 
	consider the upper right triangle matrix of $D$ leaving out the diagonal,
	$$
		\left(
		\begin{array}{cccc}
			d_{0,1} & d_{0,2} &\;\,\dots & d_{0,n-1} \\ 
			& d_{1,2}& \;\,\dots & d_{1,n-1} \\ 
			&  & \;\,\ddots & \vdots\\
			&  &  & \d_{n-2,n-1}
		\end{array}
		\right),
	$$ 
	and read the entries column by  column as a word 
	$$w(D)=(d_{0,1},d_{0,2},d_{1,2},\dots\dots,d_{0,n-1},\dots,\d_{n-2,n-1})\,.$$
	With a lexicographical ordering on the words $w(D)$ we define 
	$$
		D_1\succeq D_2\quad\Longleftrightarrow\quad w(D_1)\succeq w(D_2)
	$$  
	for distance matrices $D_1$, $D_2$. We call a distance matrix $D=(d_{i,j})_{0\le i,j<n}$ \textbf{canonical} if 
	$$
		D\succeq (d_{\tau(i),\tau(j)})\quad\forall \tau\in S_n\,.
	$$
	By $\downarrow\!\! D$ we denote the distance matrix consisting of the first $n-1$ rows and columns of $D$. 
	With this we call a distance matrix $D$ \textbf{semi-canonical} if
	$$
		\downarrow\!\! D\succeq \downarrow\!\!(d_{\tau(i),\tau(j)})\quad\forall \tau\in S_n\,.
	$$
	A canonical distance matrix is also semi-canonical. It is left to the reader to prove that each semi-canonical distance 
	matrix $D$ can be obtained by combining a canonical distance matrix $D_1$ and a semi-canonical distance matrix $D_2$, 
	see Figure \ref{fig2}. Only the distance $d_{n-1,n-2}$ is not determined by the distances of $D_1$ and $D_2$. Here we consider 
	two cases. If we combine two $(\di'-1)$-dimensional simplices to get an $\di'$-dimensional simplex Theorem \ref{theorem_menger}
	yields a biquadratic inequality for $d_{n-1,n-2}$. In the other case we can determine one or for $n=\di+2$ at most two 
	different coordinate representations of the $n$ points similar to the proof of Lemma \ref{lemma_coordinates}, calculate  
	$d_{n-1,n-2}$, and check whether it is integral. We denote the sub routine doing this by $combine(D_1,D_2)$. At first we provide an 
	algorithm to generate $\di$-dimensional integral simplices. Therefore we assume that for a
	given diameter $\Delta$, this is the largest distance, we have two lists $\mathcal{L}_{\di}^c$, $\mathcal{L}_{\di}^s$ of the 
	canonical and the semi-canonical $(\di-1)$-dimensional integral simplices with diameter $\Delta$ which are ordered 
	by $\prec$, respectively. The following algorithm determines the lists $\mathcal{L}_{\di+1}^c$ and $\mathcal{L}_{\di+1}^s$ 
	of the $\di$-dimensional integral simplices with diameter $\Delta$ ordered by $\prec$.
	\begin{algorithmus}{}\\
		\label{algo_simplices}
  		{\em Input:} $\mathcal{L}_{\di}^c$, $\mathcal{L}_{\di}^s$\\
		{\em Output:} $\mathcal{L}_{\di+1}^c$, $\mathcal{L}_{\di+1}^s$\\
		\textbf{begin}\\
		\T$\mathcal{L}_{\di+1}^c=\emptyset,\quad\mathcal{L}_{\di+1}^s=\emptyset$\\
		\T\textbf{loop over} $x\in \mathcal{L}_{\di}^c$ \textbf{do}\\
		\T\T\textbf{loop over} $\mathcal{L}_{\di}^s\ni\, y\preceq x$ \textbf{with} $\downarrow\!\! x=\downarrow\!\! y$ \textbf{do}\\
		\T\T\T\textbf{loop over} $z\in combine(x,y)$ \textbf{do}\\
		\T\T\T\T\textbf{if} $z$ is canonical \textbf{then} $\mathcal{L}_{\di+1}^c\longleftarrow z$ \textbf{end}\\
		\T\T\T\T\textbf{if} $z$ is semi-canonical \textbf{then} $\mathcal{L}_{\di+1}^s\longleftarrow z$ \textbf{end}\\
		\T\T\T\textbf{end}\\		
		\T\T\textbf{end}\\		
		\T\textbf{end}\\
		\textbf{end}\\
	\end{algorithmus}
	Because an $\di$-dimensional simplex is an $\di$-dimensional point set consisting of $n=\di+1$ points we can use Algorithm 
	\ref{algo_simplices} to generate complete lists $\mathcal{M}_{\di+1}^c$, $\mathcal{M}_{\di+1}^s$ of the canonical 
	and semi-canonical $\di$-dimensional integral point sets with diameter $\Delta$ consisting of $\di+1$ points, respectively. An 
	$\di$-dimensional point set is in semi-general position if no $\di+1$ points are situated on an $(\di-1)$-dimensional hyperplane. 
	Using Theorem \ref{characteristic} we can give an algorithm to determine the lists $\mathcal{M}_{n}^c$ and 
	$\mathcal{M}_{n}^s$ of the $\di$-dimensional integral point sets in semi-general position consisting of $n$ points 
	with diameter $\Delta$.
	\begin{algorithmus}{}\\
		\label{algo_pointsets}
		{\em Input:} $\mathcal{M}_{n-1}^c$, $\mathcal{M}_{n-1}^s$\\
		{\em Output:} $\mathcal{M}_{n}^c$, $\mathcal{M}_{n}^s$\\
		\textbf{begin}\\
		\T$\mathcal{M}_{n}^c=\emptyset,\quad\mathcal{M}_{n}^s=\emptyset$\\
		\T\textbf{loop over} $x\in \mathcal{M}_{n-1}^c$ \textbf{do}\\
		\T\T\textbf{loop over} $\mathcal{M}_{n-1}^s\ni\, y\preceq x$ \textbf{with} $\downarrow\!\! x=\downarrow\!\! y$ 
		\textbf{and} $\car(x)=\car(y)$ \textbf{do}\\
		\T\T\T\textbf{loop over} $z\in combine(x,y)$ \textbf{do}\\
		\T\T\T\T\textbf{if} $z$ is canonical \textbf{then} $\mathcal{M}_{n}^c\longleftarrow z$ \textbf{end}\\
		\T\T\T\T\textbf{if} $z$ is semi-canonical \textbf{then} $\mathcal{M}_{n}^s\longleftarrow z$ \textbf{end}\\
		\T\T\T\textbf{end}\\		
		\T\T\textbf{end}\\		
		\T\textbf{end}\\
		\textbf{end}\\
	\end{algorithmus}
  
	
	\section{Improvements}
	
	\begin{table}[!ht]
		\begin{center}
			\begin{tabular}{|r|r|r|r||r|r|r|r|}
				\hline
				$\Delta$	& $\hat{\Psi}(3,\Delta)$ & $\Psi(3,\Delta)$ & $\tilde{\alpha}(3,\Delta)$ & 
				$\Delta$	& $\hat{\Psi}(3,\Delta)$ & $\Psi(3,\Delta)$ & $\tilde{\alpha}(3,\Delta)$ \\
				\hline
   			 1 &         1 &      1 &      1 & 26 &   521610123 &   521589 &  356333 \\
				 2 &        13 &      9 &      6 & 27 &   700065646 &   629939 &  428030 \\
				 3 &       111 &     35 &     24 & 28 &   929489332 &   753113 &  510829 \\
				 4 &       602 &    149 &     70 & 29 &  1222613496 &   832969 &  605970 \\
				 5 &      2592 &    305 &    176 & 30 &  1592477593 &  1038224 &  714505 \\
				 6 &      8833 &    770 &    380 & 31 &  2059062666 &  1145517 &  838646 \\
				 7 &     26564 &   1379 &    754 & 32 &  2638060710 &  1439990 &  978820 \\
				 8 &     68800 &   2761 &   1368 & 33 &  3357319548 &  1568195 & 1137638 \\
				 9 &    162330 &   4182 &   2333 & 34 &  4241882219 &  1804079 & 1316239 \\
				10 &    353100 &   6660 &   3786 & 35 &  5323350205 &  2062374 & 1516567 \\
				11 &    719688 &  10254 &   5894 & 36 &  6638917601 &  2475320 & 1740591 \\
				12 &   1378977 &  16714 &   8839 & 37 &  8232016014 &  2613730 & 1990484 \\
				13 &   2526059 &  21902 &  12891 & 38 & 10148934902 &  3037708 & 2268149 \\
				14 &   4434103 &  30115 &  18289 & 39 & 12445587259 &  3430131 & 2575954 \\
				15 &   7490297 &  41250 &  25339 & 40 & 15183055989 &  4015829 & 2916089 \\
				16 &  12256818 &  59995 &  34436 & 41 & 18437914417 &  4224348 & 3291649 \\
				17 &  19551329 &  72315 &  46054 & 42 & 22280569281 &  4966748 & 3704516 \\
				18 &  30264028 &  96502 &  60474 & 43 & 26818516374 &  5278577 & 4158686 \\
				19 &  45952871 & 119896 &  78406 & 44 & 32132601503 &  6213243 & 4655277 \\
				20 &  68191989 & 162600 & 100277 & 45 & 38348410933 &  6821671 & 5198318 \\
				21 &  99420707 & 196490 & 126838 & 46 & 45598443859 &  7428904 & 5791458 \\
				22 & 142558111 & 245591 & 158772 & 47 & 54019488362 &  8057637 & 6437526 \\
				23 & 201289670 & 289672 & 196799 & 48 & 63756807373 &  9675353 & 7139157 \\
				24 & 279728968 & 388051 & 241672 & 49 & 75019979427 & 10055859 & 7901871 \\
				25 & 384663513 & 440140 & 294681 & 50 & 87968187078 & 11262298 & 8727553 \\
				\hline 
			\end{tabular}\\[1mm] 
			\caption{Number of calls of $combine(x,y)$.}
			\label{table_comp}
		\end{center}
	\end{table}
	
	To demonstrate the significance of Theorem \ref{characteristic} for an efficient enumeration algorithm for integral 
	point sets we compare in Table \ref{table_comp} the number $\Psi(3,\Delta)$ of calls of $combine(x,y)$ 
	in Algorithm \ref{algo_pointsets} for $\di=3$ and $n=5$ to the number $\hat{\Psi}(3,\Delta)$ of calls of 
	$combine(x,y)$ without using Theorem \ref{characteristic}. Additionally we give the number 
	$\tilde{\alpha}(3,\Delta)$ of semi-canonical integral tetrahedrons with diameter $\Delta$.
	
	
	\section{Minimum diameters}
  
	From the combinatorial point of view there is a natural interest in the minimum diameter $d(\di,n)$ of
	$\di$-dimensional integral point sets consisting of $n$ points. By $\overline{d}(\di,n)$ we denote 
	the minimum diameter of $\di$-dimensional integral point sets in semi-general position. If additionally 
	no $\di+2$ points lie on an $\di$-dimensional sphere we denote the corresponding minimum diameter by 
	$\dot{d}(\di,n)$ and say the points are in general position. To check semi-general position we can use the 
	Cayley-Menger matrix and test whether $V_{\di}=0$ or not. In the case of general position we have the following theorem.
  
	\begin{satz}{}
		Given $\di+2$ points in $\mathbb{E}^{\di}$, with pairwise distances $d_{i,j}$ and no $\di+1$ points in 
		an $\di-1$-dimensional plane, lie on an $\di$-dimensional sphere if and only if
		$$
			\left|
			\begin{array}{cccc}
				0 & d_{0,1}^2 & \dots & d_{0,\di+1}^2 \\ 
				d_{1,0}^2 & 0 & \ddots & \vdots \\ 
				\vdots & \ddots & \ddots & d_{\di,\di+1}^2 \\ 
				d_{\di+1,0}^2 & \dots & d_{\di+1,\di}^2 & 0
			\end{array} 
			\right|=0\,.
		$$
	\end{satz}
	See \cite{odd_distances,dipl_piepmeyer} for a proof.
  
	We have implemented Algorithm \ref{algo_simplices} and Algorithm \ref{algo_pointsets} and received the following 
	values for minimum diameters, see also \cite{integral_distances_in_point_sets,kurz_wassermann,dipl_piepmeyer}. The 
	values not previously known in the literature are emphasised.
	$$
		\overline{d}(3,n)_{4\le n\le 7}=\dot{d}(3,n)_{4\le n\le 7}=1,3,16,\textbf{44}.
	$$
	$$
		\overline{d}(4,n)_{5\le n\le 8}=1,\textbf{4},\textbf{11},\textbf{14}.
	$$
	$$
		\dot{d}(4,n)_{5\le n\le 8}=1,\textbf{4},\textbf{7},\textbf{14}.
	$$
	$$
		\overline{d}(5,n)_{6\le n\le 9}=\dot{d}(5,n)_{6\le n\le 9}=1,\textbf{4},\textbf{5},\textbf{8}.
	$$
	
	To determine $d(\di,n)$ we have to modify Algorithm \ref{algo_pointsets} because not every $\di+1$ points of an 
	$\di$-dimensional pointset span an $\di$-dimensional simplex. So we have to combine lower dimensional point sets with 
	$\di$-dimensional point sets. We leave the details to the reader and give only the results,
	
	\begin{eqnarray*}
		\quad\quad\quad\quad d(3,n)_{4\le n\le 23} & = &
		 1,3,4,8,13,\textbf{16},17,\textbf{17},\textbf{17},\textbf{56},\textbf{65},\textbf{77},\\
		& & \textbf{86},\textbf{99},\textbf{112}, \textbf{133}, \textbf{154}, \textbf{195}, \textbf{212}, \textbf{228}.
	\end{eqnarray*}

	\bibliography{characteristic}
  	\bibdata{characteristic}
	\bibliographystyle{plain}  
\end{document}